\documentclass[12pt,reqno]{amsart}
\usepackage{amsfonts,amsmath,amssymb}

\allowdisplaybreaks
\advance\textheight by \topskip

\usepackage{mathrsfs}

\usepackage[latin1]{inputenc}
\usepackage{graphicx}

\usepackage{tikz}
\newtheorem{theorem}{Theorem}
\usepackage{euscript}

\def\[{[\! [}
\def\]{]\! ]}

\title[ Lattice paths with infinitely many down steps]{Lattice paths with infinitely many down steps -- the negative boundary model}

\author[H.~Prodinger]{Helmut Prodinger}

\address{Helmut Prodinger,
Mathematics Department, Stellenbosch University,
7602 Stellenbosch, South Africa.}
\email{hproding@sun.ac.za}

\date{\today}

\begin{document}

\begin{abstract}
	We consider a variation of Dyck paths, where additionally to steps $(1,1)$ and $(1,-1)$  down-steps  $(1,-j)$, for $j\ge2$ are  allowed.
	We give credits to Emeric Deutsch for that.
	The enumeration of such objects living in a strip is performed. Methods are  the   kernel method and techniques from linear algebra.
\end{abstract}

\maketitle

\section{Introduction}

Emeric Deutsch \cite{Deutsch} had the idea to consider a variation of ordinary Dyck paths, by augmenting the usual up-steps and down-steps by one unit each,
by down-steps of size $3,5,7,\dots$. This leads to ternary equations, as can be seen for instance from \cite{Deutsch-ternary}.

The present author started to investigate a related but simpler model of down-steps $1,2,3,4,\dots$ and investigated it (named Deutsch paths in honour of Emeric Deutsch)
in a series of papers, \cite{Deutsch1,Deutsch-slice,Prodinger-fibo}.

This paper is an further member of this series: The condition that (as with Dyck paths) the paths cannot enter negative territory, is relaxed, by introducing
a negative boundary $-t$. Here are two recent publications about such a negative boundary: \cite{Selkirk-master} and \cite{Prodinger-negative}.

Instead of allowing negative altitudes, we think about the whole system shifted up by $t$ units, and start at the point $(0,t)$ instead. This
is much better for the generating functions that we are going to investigate. Eventually, the results can be re-interpreted as results about 
enumerations with respect to a negative boundary.

The setting with flexible initial level $t$ and final level $j$ allows us to consider the Deutsch paths also from left to right (they are not symmetric!),
without any new computations.

The next sections achieves this, using the celebrated kernel-method,  one of the tools that is dear to our heart \cite{Prodinger-kernel}.

In the following section, an additional upper bound is introduced, so that the Deutsch paths live now in a strip. The way to attack this
is linear algebra. Once everything has been computated, one can relax the conditions and let lower/upper boundary go to $\mp\infty$.

\section{Generating functions and the kernel method}

As discussed, we consider Deutsch paths starting at $(0,t)$ and ending at $(n,j)$, for $n,t,j\ge0$. First we consider univariate generating functions
$f_j(z)$, where $z^n$ stays for $n$ steps done, and $j$ is the final destination. The recursion is immediate:
\begin{equation*}
f_j(z)=[\![t=j]\!]+zf_{j-1}(z)+z\sum_{k>j}f_k(z),
\end{equation*}
where $f_{-1}(z)=0$. Next, we consider
\begin{equation*}
F(z,u):=\sum_{j\ge0}f_j(z)u^j,
\end{equation*}
and get
\begin{align*}
F(z,u)&=u^t+zuF(z,u)+z\sum_{j\ge0}u^j\sum_{k>j}f_k(z)\\
&=u^t+zuF(z,u)+z\sum_{k>0}f_k(z)\sum_{0\le j<k}u^j\\
&=u^t+zuF(z,u)+z\sum_{k\ge0}f_k(z)\frac{1-u^k}{1-u}\\
&=u^t+zuF(z,u)+\frac z{1-u}[ F(z,1)-F(z,u)]\\
&=\frac{u^t(1-u)+zF(z,1)}{z-zu+zu^2+1-u}.
\end{align*}
Since the critical value is around $u=1$, we write the denominator as
\begin{equation*}
z(u-1)^2+(u-1)(z-1)+z=z(u-1-r_1)(u-1-r_2),
\end{equation*}
 with
\begin{align*}
r_1=\frac{1-z+\sqrt{1-2z-3z^2}}{2z},\quad r_2=\frac{1-z-\sqrt{1-2z-3z^2}}{2z}.
\end{align*}
The factor $(u-1-r_2)$ is bad, so the numerator must vanish for $[u^t(1-u)+zF(z,1)]|_{u=1+r_2}$, therefore
\begin{equation*}
zF(z,1)=(1+r_2)^tr_2.
\end{equation*} 
Furthermore
\begin{equation*}
F(z,u)=
\frac{\frac{u^t(1-u)+zF(z,1)}{u-r_2}}{z(u-r_1)}.
\end{equation*}
The expressions become prettier using the substitution $z=\frac{v}{1+v+v^2}$; then
\begin{equation*}
r_1=\frac{1}{v},\quad r_2=v.
\end{equation*}
It can be proved by induction (or computer algebra) that 
\begin{equation*}
\frac{u^t(1-u)+v(1+v)^t}{u-1-v}=-v\sum_{k=0}^{t-1}(1+v)^{t-1-k}-u^t.
\end{equation*}
Furthermore
\begin{equation*}
\frac1{z(u-1-r_1)}=-\frac1{z(1+r_1)(1-\frac{u}{1+r_1})},
\end{equation*}
and so 
\begin{equation*}
f_j(z)=[u^j]F(z,u)=[u^j]\biggl[v\sum_{k=0}^{t-1}(1+v)^{t-1-k}u^k+u^t\biggr]\sum_{\ell\ge0}\frac{u^{\ell}}{z(1+r_1)^{\ell+1}}.
\end{equation*}
Of interest are two special cases: 
The case that was studied before \cite{Deutsch1} is $t=0$:
\begin{equation*}
	f_j=\frac{(1+v+v^2)v^{j}}{(1+v)^{j+1}}.
\end{equation*}
The other special case is $j=0$ for general $t$, as it may be interpreted as Deutsch paths read from right to left, starting at level $0$ and 
ending at level $t\ge1$ (for $t=0$, the previous formula applies):
\begin{align*}
	f_0(z)&=[u^0]\biggl[v\sum_{k=0}^{t-1}(1+v)^{t-1-k}u^k+u^t\biggr]\sum_{\ell\ge0}\frac{u^{\ell}}{z(1+r_1)^{\ell+1}}\\
	&=v(1+v)^{t-1}\frac{1}{z(1+r_1)}=v(1+v+v^2)(1+v)^{t-2}.
\end{align*}

The next section will present a simplification of the expression for $f_j(z)$, which could be obtained directly by distinguishing cases and summing some geometric series.

\section{Refined analysis: lower and upper boundary}

Now we consider Deutsch paths bounded from below by zero and bounded from above by $m-1$; they start at level $t$ and end at level $j$ after $n$ steps.
For that, we use generating functions $\varphi_j(z)$ (the quantity $t$ is a silent parameter here). The recursions that are straight-forwarded are best organized in a
matrix, as the following example shows.
\begin{equation*}
	\left(\begin{matrix}
		1&-z&-z&-z&-z&-z&-z&-z\\
		-z&	1&-z&-z&-z&-z&-z&-z\\
		0&	-z&	1&-z&-z&-z&-z&-z\\
		0& 0&	-z&	1&-z&-z&-z&-z\\
		0& 0& 0&	-z&	1&-z&-z&-z\\
		0& 0& 0&0&	-z&	1&-z&-z\\
		0& 0& 0&0&0&	-z&	1&-z\\
				0& 0& 0&0&0&0&	-z&	1\\
	\end{matrix}\right)
	\left(\begin{matrix}
		\varphi_0\\
		\varphi_1\\
		\varphi_2\\
		\varphi_3\\
		\varphi_4\\
				\varphi_5\\
						\varphi_6\\
												\varphi_7\\
	\end{matrix}\right)=
	\left(\begin{matrix}
		0\\
		0\\
		0\\
		1\\
		0\\
			0\\
				0\\
				0\\
	\end{matrix}\right)
\begin{tikzpicture} 
\draw [](0,0)--(0,0);
\node at (0.5,1.3){\Bigg\}$t$};
\end{tikzpicture}
\end{equation*}	

The goal is now to solve this system. For that the substitution $z=\frac{v}{1+v+v^2}$ is used throughout. The method is to use Cramer's rule, which means that the right-hand side has to replace various columns of the matrix, and determinants have to be computed. At the end, one has to divide by the determinant of the system.

Let $D_m$ be the determinant of the matrix with $m$ rows and columns. 
The recursion
\begin{equation*}
	(1+v+v^2)^2m_{n+2}-(1+v+v^2)(1+v)^2D_{m+1}+v(1+v)^2D_{m}=0
\end{equation*}	
appeared already in \cite{Deutsch1} and is not difficult to derive and to solve:
\begin{equation*}
	D_m=\frac{(1+v)^{m-1}}{(1+v+v^2)^m}\frac{1-v^{m+2}}{1-v}.
\end{equation*}

To solve the system with Cramer's rule, we must compute a determinant of the following type,
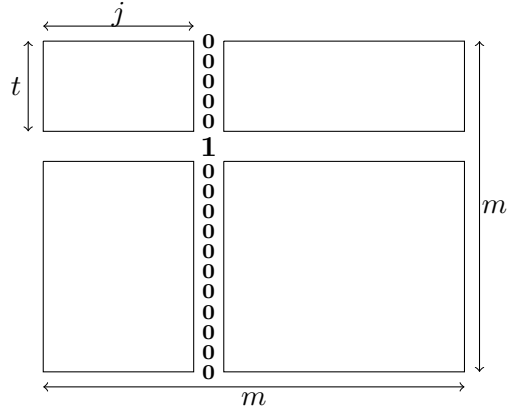
\begin{center}\small
	\begin{tikzpicture}
		[scale=0.4]
		\draw (0,0)--(5,0)--(5,-3)--(0,-3)--(0,0);
		\node at (2.5,0.9){$j$};
		\draw[<->](0,0.5)--(5,0.5);
		\draw[<->](-0.5,0)--(-0.5,-3);
		\node at (-0.9,-1.5){$t$};
		\node at (5.5,-3.5){$\boldsymbol{1}$};	
		\newcommand\x{6};\draw (0+\x,0)--(8+\x,0)--(8+\x,-3)--(0+\x,-3)--(0+\x,0);	
\newcommand\y{-4};				\draw (0,0+\y)--(5,0+\y)--(5,-7+\y)--(0,-7+\y)--(0,0+\y);
			\draw (0+\x,0+\y)--(8+\x,0+\y)--(8+\x,-7+\y)--(0+\x,-7+\y)--(0+\x,0+\y);
			\draw[<->](0,-11.5)--(14,-11.5);
						\draw[<->](14.5,0)--(14.5,-11);
			\node at (7,-11.9){$m$};
						\node at (15.0,-5.5){$m$};
			\foreach \x in {0,1,2,3,4}
			{
\node at (5.5,-\x/1.5){$\tiny\boldsymbol{0}$};
			}

		\foreach \x in {6.5,7.5,8.5,9.5,10.5,11.5,12.5,13.5,14.5,15.5,16.5}
		{
			\node at (5.5,-\x/1.5){$\tiny\boldsymbol{0}$};
		}
			
	\end{tikzpicture}
\end{center}
where the various rows are replaced by the right-hand side. While it is not impossible to solve this recursion by hand,
it is very easy to make mistakes, so it is best to employ a computer. Let $D(m;t,j)$ the determinant according to the
drawing.

It is not unexpected that the results are different for $j<t$ resp.\ $j\ge t$. Here is what we found:
\begin{equation*}
D(m;t,j)=\frac{(1+v)^{t-j-3+m}(1-v^{j+1})v(1-v^{m-t})}{(1-v)^2(1+v+v^2)^{m-1}},\quad\text{for}
\ j<t,
\end{equation*}
\begin{equation*}
	D(m;t,j)=\frac{v^{j-t}(1-v^{t+2})(1-v^{1-j+m})}{(1-v)^2(1+v+v^2)^{m-1}(1+v)^{j-t+3-m}},\quad\text{for}
	\ j\ge t.
\end{equation*}
To solve the system, we have to divide by the determinant $D_m$, with the result
\begin{equation*}
\varphi_j=	\frac{D(m;t,j)}{D_m}=
\frac{(1+v)^{t-j-2}(1-v^{j+1})v(1-v^{m-t})(1+v+v^2)}{(1-v)(1-v^{m+2})}
,\quad\text{for}
	\ j<t,
\end{equation*}
\begin{equation*}
\varphi_j=	\frac{D(m;t,j)}{D_m}=\frac{v^{j-t}(1-v^{t+2})(1-v^{1-j+m})(1+v+v^2)}{(1-v)(1+v)^{j-t+2}(1-v^{m+2})}
,\quad\text{for}
	\ j\ge t.
\end{equation*}
We found all this using Computer algebra. Some critical minds may argue that this is only experimental. One way of rectifying this
would be to show that indeed the functions $\varphi_j$ solve the system, which consists of summing various geometric series; again,
a computer could be helpful for such an enterprise.

Of interest are also the limits for $m\to\infty$, i.e., no upper boundary:
\begin{equation*}
	\varphi_j=\lim_{m\to\infty}\frac{D(m;t,j)}{D_m}=
	\frac{(1+v)^{t-j-2}(1-v^{j+1})v(1+v+v^2)}{(1-v)}
	,\quad\text{for}
	\ j<t,
\end{equation*}
\begin{equation*}
	\varphi_j=\frac{v^{j-t}(1-v^{t+2})(1+v+v^2)}{(1-v)(1+v)^{j-t+2}	}
	,\quad\text{for}
	\ j\ge t.
\end{equation*}
The special case $t=0$ appeared already in the previous section:
\begin{equation*}
	\varphi_j=\frac{v^{j}(1+v+v^2)}{(1+v)^{j+1}	}.
\end{equation*}
Likewise, for $t\ge1$,
\begin{equation*}
	\varphi_0=v(1+v+v^2)(1+v)^{t-2}.
\end{equation*}
In particular, the formul\ae\ show that the expression from the previous section can be simplified in general, which
could have been seen directly, of course.

\begin{theorem}
	The generating function of Deutsch path with lower boundary 0 and upper boundary $m-1$, starting at $(0,t)$ and ending at $(n,j)$ is given by
	\begin{gather*}
\frac{(1+v)^{t-j-2}(1-v^{j+1})v(1-v^{m-t})(1+v+v^2)}{(1-v)(1-v^{m+2})}
,\quad\text{for}
\ j<t,\\
\frac{v^{j-t}(1-v^{t+2})(1-v^{1-j+m})(1+v+v^2)}{(1-v)(1+v)^{j-t+2}(1-v^{m+2})}
,\quad\text{for}
\ j\ge t,		
	\end{gather*}
with the substitution $z=\dfrac{v}{1+v+v^2}$.
\end{theorem}
By shifting everything down, we can interpret the results as Deutsch walks between boundaries $-t$ and $m-1-t$, starting at the origin $(0,0)$ and ending at $(n,j-t)$.

\begin{theorem}
	The generating function of Deutsch path with lower boundary $-t$ and upper boundary $h$, starting at $(0,0)$ and ending at $(n,i)$ with $-t\le i\le h$ is given by
	\begin{gather*}
		\frac{(1+v)^{i-2}(1-v^{i+t+1})v(1-v^{h+1})(1+v+v^2)}{(1-v)(1-v^{h+t+3})}
		,\quad\text{for}
		\ i<0,\\	
		\frac{v^{i}(1-v^{t+2})(1-v^{2-i+h})(1+v+v^2)}{(1-v)(1+v)^{i+2}(1-v^{h+t+3})}
		,\quad\text{for}
		\ i\ge 0.		
			\end{gather*}
\end{theorem}
It is possible to consider the limits $t\to\infty$ and/or $h\to\infty$ resulting in simplified formul\ae.

\section{Conclusion}

Various parameters could be worked out starting from the present findings. Currently, nothing to that effect has been done.


\bibliographystyle{plain}


\end{document}